\documentstyle[12pt]{article}

\newfont{\num}{msym10}
\newfont{\g}{eufm10}

\newcommand{\gtm}{\mbox{\g m}}

\newcommand{\nc}{\mbox{\num C}}
\newcommand{\nr}{\mbox{\num R}}
\newcommand{\nz}{\mbox{\num Z}}
\newcommand{\none}{\mbox{\bf 1}}

\newtheorem{lm}{Lemma}
\newtheorem{prop}{Proposition}
\newtheorem{th}{Theorem}

\title{Complex Divisors on Algebraic Curves and Some Applications to String
Theory}
\author{A.~A.~Voronov\thanks{This paper is in final form and no version of it
will be submitted for publication elsewhere.} \\
Department of Mathematics, Moscow State University, \\ Moscow, 117234,
USSR\thanks{Present address: Department of Mathematics, Princeton University,
Princeton, NJ
08544, USA.}}

\date{February 26, 1992}

\begin{document}
\bibliographystyle{plain}
\maketitle

This talk presents some new notions of the theory of complex algebraic curves
which have appeared as algebraic tools in string theory (see \cite{me:unif} for
more details). In a sense, we have materialized non-existed complex powers of
invertible sheaves on algebraic curves introduced at the level of the Atiyah
algebras of invertible sheaves by Beilinson and Schechtman \cite{bs}. The {\em
Atiyah algebra\/} $A_{\cal L}$ in the case of an invertible sheaf $\cal L$ over
a complete complex algebraic curve $X$ is just the sheaf of differential
operators of order $\leq 1$ on $\cal L$. There takes
place the exact sequence $0 \rightarrow {\cal O} \rightarrow A_{\cal L}
\rightarrow T \rightarrow 0$, where $\cal O$ is the structural sheaf and $T$ is
the tangent sheaf of $X$. The diagram
\[\begin{array}{ccccccc}
0 \rightarrow & {\cal O} & \rightarrow & A_{\cal L} & \rightarrow & T &
\rightarrow 0 \\
 & c \downarrow & & \downarrow & & \parallel \\
0 \rightarrow & {\cal O} & \rightarrow & ? & \rightarrow & T & \rightarrow 0
\end{array},\]
where $c$ is the operator of multiplication by $c \in \nc$, as usual, can be
completed to commutative with a sheaf $A_{{\cal L}^c}$ interpreted as the
Atiyah algebra of a (nonexisting) invertible sheaf ${\cal L}^c$, $\cal L$ to
the power $c$. The Atiyah algebra keeps incomplete information about its
invertible sheaf, because an isomorphism $A_{\cal L}
\stackrel{\sim}{\rightarrow} A_{\cal O}$ implies the existence of a canonical
flat holomorphic connection on $\cal L$, but not the triviality of $\cal L$. It
is remarkable that using Atiyah algebras one can apply local arguments to
Riemann-Roch type theorems.

Here we deal with really existing invertible sheaves, corresponding to divisors
with complex coefficients, but having integral degree. In particular, every
invertible sheaf of degree 0 can be raised to a complex power.

\section{Complex Divisors}

Let $X$ be a complete complex curve of genus $g$, with a fixed ordered set
$\gtm = \{Q_1, \ldots, Q_n\}$ of $n$ distinct points on $X$ and a closed disk
$B$, considered up to an isotopy in $X \setminus \gtm$, such that $\gtm \subset
B$. A {\em complex divisor\/}  is a formal sum
\[  D = \sum_{P \in X} n_P \cdot P, \]
where
\[ n_P \in \left\{ \begin{array}{ll}
                                     \nc  & \mbox{for $P \in \gtm$}, \\
                                     \nz  & \mbox{otherwise},
                                   \end{array}
                           \right.\]
\[   \deg D := \sum_{P \in X} n_P \; \in \nz, \]
and only a finite number of $n_P \neq 0$. The corresponding {\em group of
complex divisors\/} is denoted by Div$\, (X,\gtm,B)$.

This definition may be not interesting itself, but it leads to a new class of
invertible sheaves over $X$.

\section{Multiple Valued Meromorphic Functions}

Let $p: \widetilde{X \setminus \gtm} \rightarrow X \setminus \gtm$ be the
universal covering with the complex structure lifted from the base. Denote by
$H$ the kernel of the natural epimorphism $\pi_1(X \setminus \gtm) \rightarrow
\pi_1(X)$ determined by the embedding  $X \setminus \gtm \hookrightarrow X$:
\[1 \rightarrow  H  \rightarrow  \pi_1(X \setminus \gtm)  \rightarrow  \pi_1(X)
 \rightarrow 1 .\]
We will call a holomorphic function $\phi$ on $\widetilde{X \setminus \gtm}$
(more precisely, a section of the sheaf $p_*{\cal O}_{\widetilde{X \setminus
\gtm}}$) a {\em multiple valued holomorpic function\/} on $X$, if
\begin{enumerate}
  \item $\phi$ is  $\pi_1(X \setminus B)$-invariant (i.e., $\phi$ is single
valued outside $B$),
  \item for every $\sigma \in H$ \  $\phi^\sigma = f_\sigma \cdot \phi$, where
$\phi^\sigma(z) := \phi(\sigma z)$, and $f_\sigma$ is a constant ($f_\sigma$ is
called the {\em multiplicator\/}),
  \item the branches of $\phi$, as the branches of a multiple valued analytic
function on $X$, have only removable singularities in \gtm, that is for any
$Q_i \in \gtm$ and any sequence $\{a_m\}$ in a domain of univalence in
$\widetilde{X \setminus \gtm}$, such that $p(a_m) \rightarrow Q_i$ as $m
\rightarrow \infty$, there exists a limit $\lim_{m \rightarrow \infty}
\phi(a_m)$, depending only on $Q_i$:
\[  \phi(Q_i) := \lim_{m \rightarrow \infty} \phi(a_m)  .  \]
\end{enumerate}
Denote by $\cal O'$ the sheaf of holomorphic multiple valued functions on $X$.
Denote the corresponding sheaf of fields of fractions by $\cal M'$. We will
call sections of $\cal M'$ {\em multiple valued meromorphic functions\/} on
$X$. The following simple lemma describes the local behavior of such functions.
\begin{lm}
Let $z$ be a holomorphic coordinate on $X$ near $Q_i \in \gtm \subset X$.
\begin{enumerate}
 \item If $\phi \in \cal O'$, then either
\[ \phi(z) = z^A \cdot \sum^\infty_{j=0} \alpha_j z^j, \mbox{ where } 0 < \Re A
\leq 1, \]
or
\[ \phi(z) = \sum^\infty_{j=0} \alpha_j z^j.  \]
 \item If $\phi \in \cal M'$, then
\[ \phi(z) = z^A \cdot \sum^\infty_{j=n_0} \alpha_j z^j, \mbox{ where } 0 \leq
\Re A < 1. \Box \]
\end{enumerate}
\end{lm}
{\em Note}. One should remember that these expansions may also get monodromy at
other points $Q_i \in \gtm$.
\smallskip \\
{\em Definition}. The number $A+n_0$ is called the order ord$_{Q_i} \phi$ of
the multiple valued holomorphic function $\phi$ at the singular point $Q_i$.

Let $\phi \in \Gamma(X, {\cal M}')$ be a globally defined multiple valued
holomorphic function on $X$. Then $\sum_{P \in X} \mbox{ord}_P \phi = 0$,
because d$\, \log \phi$ is a differential of the third kind on $X$ and the sum
of its residues vanishes.
\smallskip\\
{\em Definition}. A divisor of the type
\[
\mbox{div}\, \phi := \sum_{P \in X} \mbox{ord}_P \phi \cdot P  \]
is called principal.
\smallskip\\
Define the {\em group\/} Cl$\,(X, \gtm, B)$  {\em of classes of complex
divisors\/} as the quotient-group of the group Div$\, (X, \gtm, B)$ by the
subgroup of principal divisors.

\section{Complex Divisors and Invertible Sheaves}

\begin{prop}
\begin{sloppypar}
\begin{enumerate}
 \item The group {\rm Cl}$\,(X, \gtm, B)$ is isomorphic to the group ${\rm
Cl}\,(X)$ of classes of ordinary (integral) divisors on $X$.
 \item The group ${\rm Div}\,(X, \gtm, B)$ is isomorphic to the group of
invertible $\cal O$-submodules in $\cal M'$.
\end{enumerate}
\end{sloppypar}
\end{prop}
{\em Proof}. Part 1 evidently follows from 2, so let us prove 2. Choose a
covering of $X$ with two open subsets $U_1 := \{ \mbox{a $\delta$-neighborhood
of $B$ for small $\delta > 0$}$, $U_2 := X \setminus B$, and given complex
divisor $D = \sum n_P \cdot P$ take a multiple valued meromorphic function
$f_1$ on $U_1$, such that ord$_P f_1 = n_P$ for $P \in U_1$, and a multiple
valued meromorphic function $f_2$ on $U_2$, such that ord$_P f_2 = n_P$ for $P
\in U_2$ and $f_2^{\gamma_i} = \exp (2 \pi \sqrt{-1} \cdot n_{Q_i}) \cdot
f_2,\; i = 1, \ldots, n$, where $\gamma_i$ is a loop in $B$ containing the
single point $Q_i$. Then $f_1 / f_2$ is a single valued nonzero holomorphic
function on $U_1 \cap U_2$, i.e., $f_1 / f_2 \in \Gamma(U_1 \cap U_2, {\cal
O}^*)$, and it determines an $\cal O$-submodule ${\cal O}(D)$ in $\cal M'$
having $f_2 / f_1$ as the glueing function. Thus, ${\cal O}(D)$ is an ordinary
invertible sheaf on $X$. $\Box$

\section{The Weil-Deligne Pairing}

\begin{sloppypar}
Let ${\cal L}_1$, ${\cal L}_2$ be two invertible $\cal O$-modules. [They may
well be $\cal O$-submodules in $\cal M'$]. Define a complex vector space
$\langle {\cal L}_1, {\cal L}_2 \rangle$ as the space generated by the
expressions
\begin{equation}
 \langle l_1, l_2 \rangle ,  \label{wd}
\end{equation}
where $l_1$ and $l_2$ are single valued (i.e., having integral divisors)
meromorphic sections of ${\cal L}_1$ and ${\cal L}_2$, respectively, with
nonintersecting divisors. We place the following relations on the symbols
(\ref{wd}):
\[ \langle f \cdot l_1, l_2 \rangle = f({\rm div}\, l_2) \cdot \langle l_1, l_2
\rangle , \]
\[ \langle  l_1, g \cdot l_2 \rangle = g({\rm div}\, l_1) \cdot \langle l_1,
l_2 \rangle , \]
where $f$ and $g$ are single valued meromorphic functions such that $f({\rm
div}\, l_2) := \prod_{P \in X} f(P)^{\mbox{\scriptsize ord}_P l_2} \neq 0,
\infty$ in the former formula and $g({\rm div}\, l_1) := \prod_{P \in X}
g(P)^{\mbox{\scriptsize ord}_P l_1} \neq 0, \infty$ in the latter. The
correctness of this definition is provided by Weil's reciprocity law:
\[ f({\rm div}\, g)= g({\rm div}\, f). \]
One can easily see that the space $\langle {\cal L}_1, {\cal L}_2 \rangle$ is a
one-dimensional complex vector space. We will call it the {\em Weil-Deligne
pairing\/} of ${\cal L}_1$ and ${\cal L}_2 $.
\end{sloppypar}

\section{The Arakelov-Deligne Metric}

Now, let ${\cal L}_1$ and ${\cal L}_2 $ be two Hermitian holomorphic line
bundles. Then one can define a natural Hermitian metric on the space $\langle
{\cal L}_1, {\cal L}_2 \rangle$ (cf.\  Deligne \cite{d:d}). That means that for
any two sinle valued sections $l_1$, $l_2$ of ${\cal L}_1$ and ${\cal L}_2 $
with nonintersecting divisors, there is defined a real number
\[    \parallel \langle l_1, l_2 \rangle \parallel  \in \nr .\]
Below we define an analogous metric in a more general case, when $l_1$ and
$l_2$ are not necessarily single valued, but of degree 0. We will use this
construction in string theory later.

The definition is
\begin{equation}
 \parallel \langle l_1, l_2 \rangle \parallel := \sqrt{ \prod_i
G^{\overline{n}_i}_{\mbox{\scriptsize div}\,l_2} (P_i) \cdot
G^{n_i}_{\overline{\mbox{\scriptsize div}\,l_2}} (P_i)},  \label{ad}
\end{equation}
where $\overline{\mbox{div}\,l_2}$ means the divisor with complex conjugated
coefficients, div$\, l_1 = \sum_i n_i P_i$ and $G_D (z) := \exp g_D(z)$,
$g_D(z)$ being the Green function of the divisor $D$, which is defined up to a
constant similar to the case of integral $D$ (for example, put $g_D (z) := \Re
\int_{z_0}^z \omega_D$, where $\omega_D$ is the differential of the third kind
associated with $D$, cf.\ Lang \cite{l}). The result does not depend on the
choice of Green function, because we assume deg$\, l_1 = \sum n_i = 0$.
Moreover, the obtained symbol $\parallel \langle \; , \; \rangle \parallel$ is
symmetric:
\[   \parallel \langle l_1, l_2 \rangle \parallel = \parallel \langle l_2, l_1
\rangle \parallel . \]
This can be observed from the formula
\begin{equation}
 \parallel \langle l_1, l_2 \rangle \parallel = \sqrt{ \prod_i
G^{\overline{n}_i}_{\mbox{\scriptsize div}\,l_2} (P_i) \cdot \prod_j
G^{\overline{n}'_j}_{\mbox{\scriptsize div}\,l_1} (P'_j)},  \label{adsym}
\end{equation}
where div$\, l_2 = \sum_j n'_j P'_j$. In fact, $G_{\mbox{\scriptsize div}\,l_1}
(z) = \prod_i G^{n_i}_{P_i} (z)$ and $G_P(Q) = G_Q(P)$, so (\ref{adsym}) is
equivalent to (\ref{ad}). These arguments also imply the formula
\begin{equation}
 \parallel \langle l_1, l_2 \rangle \parallel =  \prod_{i,j} G^{\Re (n_i
\overline{n}'_j)}_{P_i} (P'_j) .  \label{ad3}
\end{equation}

If $l_1$, $l_2$ and $k$ are sections of Hermitian line bundles ${\cal L}_1$ and
${\cal L}_2 $ and $\cal K$, then
\[ \parallel \langle l_1 \otimes l_2, k \rangle \parallel = \parallel \langle
l_1, k \rangle \parallel \cdot \parallel \langle l_2, k \rangle \parallel  \]
whenever both sides are defined. There are some special properties of complex
divisors: if supp$\, D_1$, supp$\, D_2 \subset \gtm$, then
\[ \parallel \langle \none_{\alpha D_1}, \none_{D_2} \rangle \parallel =
\parallel \langle \none_{D_1}, \none_{D_2} \rangle \parallel ^\alpha \mbox{ for
} \alpha \in \nr \]
and
\[ \parallel \langle \none_{\alpha D_1}, \none_{D_2} \rangle \parallel =
\parallel \langle \none_{D_1}, \none_{\overline{\alpha} D_2} \rangle \parallel
\mbox{ for } \alpha \in \nc. \]

Thereby, the symbol $\parallel \langle \; , \; \rangle \parallel$ is Hermitian.
More precisely, it is the modulus of the exponent of a Hermitian form on the
vector space of complex divisors of degree 0 with support in \gtm. This
Hermitian form is easy to write out (cf.\  (\ref{ad3})):
\[ \sum_{i,j} n_i \, \overline{n}'_j \, g_{Q_i}(Q_j)  .   \]

\section{The Deligne-Riemann-Roch Theorem}

Let us consider an algebraic family of objects $(X, \gtm, B)$ over a base $S$,
i.e., a smooth projective morphism $\pi: X \rightarrow S$ of smooth complex
algebraic varieties with fiber being a connected complex curve, \gtm\ being the
disjoint union of $n$ regular sections of $\pi$ and $B$ varying continuously
along $S$. Let $D$ and $D'$ be two families of complex divisors on $X
\rightarrow S$, more generally, two invertible $\cal O$-submodules ${\cal L}_1$
and ${\cal L}_2$ in $\cal M'$. Suppose they are metrized as well as the sheaf
$\Omega$ of relative 1-differentials along the fibers of $\pi$. Then the
sheaves $\det \nr \pi_* \cal L$ (the determinant sheaf) and $\langle {\cal
L}_1, {\cal L}_2 \rangle$ (the Weil-Deligne sheaf, whose fiber over a single
curve $X_s$, $s \in S$, in the family is defined in Section~4) are defined.
$\det \nr \pi_* \cal L$ can be endowed with a Hermitian metric according to
Quillen (see \cite{d:d}), and $\langle {\cal L}_1, {\cal L}_2 \rangle$ is
metrized in Section~5. The following theorem is important for our string
applications.
\begin{th}[Deligne \cite{d:d}]
There is a canonical isometry
\[ \det \nr \pi_* ({\cal L})^2 \otimes \det \nr \pi_* ({\cal O})^{-2} =
\langle {\cal L} \otimes \Omega^* , {\cal L} \rangle  .  \Box \]
\end{th}

\section{String Applications}

The $g$-loop contribution to the string partition function can be reduced to
the integral
\[ Z_g := \int_{{\cal M}_g} \mbox{d}\pi_g  \]
of the {\em Polyakov measure\/} d$\pi_g$ over the moduli space $ {\cal M}_g$ of
complete complex algebraic curves of genus $g$. The {\em Belavin-Knizhnik
theorem\/} represents d$\pi_g$ as the modulus squared
\[ {\rm d}\pi_g = \mu_g \wedge \overline{\mu}_g   \]
of a {\em Mumford form\/} $\mu_g$, which is a section of the sheaf $\lambda_2
\otimes \lambda_1^{-13}$, where $\lambda_i := \det \nr \pi_* (\Omega^{\otimes
i})$, $\pi$ being the universal curve $\pi: X \rightarrow {\cal M}_g$.

The {\em tachyon scattering amplitude\/} is the integral
\[  A(g; {\bf p}_1, \ldots, {\bf p}_n) := \int_{{\cal M}_{g,n}}
\mbox{d}\pi_{g,n},  \]
where ${\cal M}_{g,n}$ is the moduli space of algebraic curves of genus $g$
with $n$ punctures and the measure $\mbox{d}\pi_{g,n}$ is expressed in terms of
determinants of Laplace operators and their Green functions. The vectors ${\bf
p}_i$ on which the amplitude depends are regarded as momentum vectors at the
scattering points, so they lie in the space-time of the critical dimension,
which we identify with $\nc^{13}$ endowed with the standard Hermitian metric.
These vectors satisfy the conditions:
\begin{enumerate}
 \item $\sum_{i=1}^n {\bf p}_i = 0$ (the momentum conservation law).
 \item The Hermitian square $({\bf p}_i, {\bf p}_i)$ is equal to 1 for every
$i$ (the mass of tachyon is $\sqrt{-1}$).
\end{enumerate}

Our application to string theory consists in proving the following analogue of
the Belavin-Knizhnik theorem for string amplitudes.
\begin{th}
\[{\rm d}\pi_{g,n} = \mu_{g,n,B} \wedge \overline{\mu}_{g,n,B} / \parallel
\mu_{g,n,B} \parallel^2,  \]
where $\mu_{g,n,B}$ is a local holomorphic section of the Hermitian line bundle
$\lambda_2 \otimes \lambda_1^{-13} \otimes (\bigotimes_{\nu=1}^{13} \langle
{\cal O}(D^\nu), {\cal O}(D^\nu) \rangle )^{-1} $ over the moduli space ${\cal
M}_{g,n,B}$ of the data $(X, Q_1, \ldots, Q_n, B)$. Here $D^\nu := \sum_{i=1}^n
p_i^\nu \cdot Q_i$ is  the complex divisor with the momentum components as
coefficients. The section $\mu_{g,n,B}$ is defined locally up to a holomorphic
factor. $\Box$
\end{th}


\begin{thebibliography}{1}

\bibitem{bs}
A.~A. Beilinson and V.~V. Schechtman.
\newblock Determinant bundles and {V}irasoro algebras.
\newblock {\em Commun.\ Math.\ Phys.}, 118:651--701, 1988.

\bibitem{d:d}
P.~Deligne.
\newblock Le determinant de la cohomologie.
\newblock {\em Contemp.\ Math.}, 67:93--178, 1987.

\bibitem{l}
S.~Lang.
\newblock {\em Fundamentals of {D}iophantine Geometry}.
\newblock Springer-Verlag, New York--Berlin, 1983.

\bibitem{me:unif}
A.~A. Voronov.
\newblock A unified approach to string scattering amplitudes.
\newblock {\em Commun.\ Math.\ Phys.}, 131:179--218, 1990.
\newblock Erratum: 140:415--416, 1991.

\end{thebibliography}
\end{document}